\documentclass[authoryear]{elsarticle}

\usepackage[utf8]{inputenc}

\usepackage{amsmath}
\usepackage{bm}
\usepackage{tikz}
\usepackage{pgfplots}
\usepackage{color, colortbl}

\usepackage{hyperref}

\newtheorem{thm}{Theorem}

\newproof{pf}{Proof}

\journal{arXiv.org}

\begin{document}

\begin{frontmatter}

\title{Iterative computational identification of a spacewise dependent the source in a parabolic equations}

\author[nsi,univ]{Petr N. Vabishchevich\corref{cor}}
\ead{vabishchevich@gmail.com}

\address[nsi]{Nuclear Safety Institute, Russian Academy of Sciences, 52, B. Tulskaya, Moscow, Russia}
\address[univ]{North-Eastern Federal University, 58, Belinskogo, Yakutsk, Russia}

\cortext[cor]{Corresponding author}

\begin{abstract}
Coefficient inverse problems related to identifying the right-hand side of an equation with use of additional information is of interest among inverse problems for partial differential equations.
When considering non-stationary problems, tasks of recovering the dependence of the right-hand side on time and spatial variables can be treated as independent. 
These tasks relate to a class of linear inverse problems, which sufficiently simplifies their study. 
This work is devoted to a finding the dependence of right-hand side of multidimensional parabolic equation on spatial variables using additional observations of the solution at the final point of time --- the final overdetermination.
More general problems are associated with some integral observation  of the solution on time --- the integral overdetermination.
The first method of numerical solution of inverse problems is based on iterative solution of boundary value problem for time derivative with non-local acceleration.
The second method is based on the known approach with iterative refinement of desired dependence of the right-hand side on spacial variables.
Capabilities of proposed methods are illustrated by numerical examples for model two-dimensional problem of identifying the right-hand side of a parabolic equation. 
The standard finite-element approximation on space is used, the time discretization is based on fully implicit two-level schemes.
\end{abstract}


\end{frontmatter}

\section{Introduction}

The mathematical modeling of many applied problems of science and engineering leads to the need for the numerical solution of inverse problems.
The inverse problems for partial differential equations are particularly noteworthy \cite{alifanov2011inverse,lavrentev1986ill}.
Inverse problems are formulated as non-classical problems for partial differential equations. 
They are often classified as ill-posed (conditionally well-posed) problems.
In the theoretical study, the  fundamental questions of uniqueness of the solution and its stability are primarily considered. 

Coefficient inverse problems related to identifying coefficient and/or the right-hand side of an equation with use of some additional information is of interest among inverse problems for partial differential equations.
When considering non-stationary problems, tasks of recovering the dependence of the right-hand side on time or spatial variables can be usually treated as independent. 
These tasks relate to a class of linear inverse problems, which sufficiently simplifies their study \cite{isakov2006inverse,prilepko2000methods}.
Only in some cases we have linear inverse problems --- identification of the right-hand side of equation, 
Other coefficient inverse problems are nonlinear, that significantly complicated their study.

The task of identifying the dependence of the right-hand side on spatial variables is one of the most important problems.
Additional conditions are often formulated using the solution at the final moment of time --- final overdetermination.
In more general case the overdetermination condition is stated as some time integral average --- integral overdetermination. 
The existence and uniqueness of the solution to such an inverse problem and well-posedness of this problem in various functional classes are examined, for example, in the works \cite{rundell1980determination,prilepko1987solvability,isakov1991inverse,prilepko1993certain,kamynin2005inverse}.

In the numerical solution of inverse problem the main focus is on the development of stable computational algorithms that take into account the peculiar properties of inverse problems \cite{vogel2002computational,samarskii2007numerical}.
Inverse problems for partial differential equations can be formulated as optimal control problems \cite{lions1971optimal,maksimov2002dynamical}.
Computational algorithms are based on using gradient iterative methods for corresponding residual functional \cite{johansson2007determination,d2012reconstruction,hasanov2014unified}. 
The implementation of such approaches relates to the solution of initial-boundary problems for the original parabolic equation and its conjugate. 

For the required right-hand side of a parabolic equation, which does not depend on time, an inverse problem with final overdetermination can be formulated as a boundary problem for evolution equation of the second order. 
In this case, we can use standard computational algorithms for the solution of stationary boundary value problems.
Such direct computational algorithms based on finite-difference approximation is described in the book \cite{samarskii2007numerical} (section 6.4).
In the work \cite{xiangtuan2011direct} the identification problem is numerically solved on the basis of transition to a evolutionary problem for the derivative of the solution with respect to time, peculiarity of which is the non-local boundary condition.

In this work, we construct special iterative methods for approximate solution of identification problem of  a spacewise dependence of the source in a parabolic equations.
They fully take into account  considered inverse problems features, which relate to their evolutionary character. These methods are based on the numerical solution of the standard Cauchy problems on each iteration.
The first method is based on the iterative refinement of initial condition for time derivative of the solution.
The second method relates to the iterative refinement of the dependence of the right-hand side on the spatial variables.
Such approach have been used before \cite{prilepko1993certain,PrilepkoKostin,Kostin2013}. 

The paper is organized as follows.
Statements of direct and inverse problems for second order parabolic equation are given in Section 2.
We consider identification problem of the right-hand side of two-dimensional parabolic equation, which does not depend on time. The additional information about the solution of the equation is given at final moment of time.
The finite-element approximation in space is used. 
The method of approximate solution of inverse problem based on the iterative solution of evolutionary problem for time derivative of the solution with non-local boundary conditions is considered in Section 3.
The computational algorithms of iterative identification of the right-hand side of parabolic equation is discussed in Section 4.
In  Section 5 we describe the algorithms for solving inverse problems with integral overdetermination. 
More general inverse problems, which relate to known dependence of the right-hand side on time, are considered in  Section 6.
Results of computational experiment for model boundary value problem are represented in Section 7.
Results of the work are summed up in Section 8.

\section{Problem formulation}

For simplicity, we restrict ourselves to a 2D problem. Generalization to the 3D case is trivial. 
Let ${\bm x} = (x_1, x_2)$ and $\Omega$ be a  bounded polygon.
The direct problem is formulated as follows.
We search $u({\bm x},t)$, $0 \leq t \leq T, \ T > 0$ such that 
it is the solution of the parabolic equation of second order:
\begin{equation}\label{eq:1}
  \frac{\partial u}{\partial t}- {\rm div} (k({\bm x}) {\rm grad}u) + c({\bm x}) u = f({\bm x}),
  \quad {\bm x} \in \Omega,
  \quad 0 < t \leq T ,   
\end{equation} 
with coefficients $0 < k_1 \leq k({\bm x}) \leq k_2$, $c({\bm x}) \geq 0$.
The boundary conditions are also specified:
\begin{equation}\label{eq:2}
  k({\bm x}) \frac{\partial u}{\partial n} + \mu({\bm x}) u = 0,
  \quad {\bm x} \in \partial\Omega,
  \quad 0 < t \leq T,    
\end{equation} 
where $\mu ({\bm x}) \geq \mu_1 > 0, \  {\bm x} \in \partial \Omega$ and $n$ is the normal to $\Omega$.
The initial conditions are
\begin{equation}\label{eq:3}
  u({\bm x}, 0) = u_0({\bm x}),
  \quad {\bm x} \in \Omega .  
\end{equation} 
The formulation (\ref{eq:1})--(\ref{eq:3}) presents the direct problem, 
where the right-hand side, coefficients of the equation as well as  the boundary and initial conditions are specified.

Let us consider the inverse problem, where in equation (\ref{eq:1}), the right-hand side
$f({\bm x})$ is unknown. An additional condition is often formulated as 
\begin{equation}\label{eq:4}
  u({\bm x}, T) =u_T ({\bm x}),
  \quad {\bm x} \in \Omega .  
\end{equation} 
In this case, we speak about the final overdetermination.

We assume that the above inverse problem of finding a pair of $u({\bm x},t), \ f({\bm x})$ from equations  
(\ref{eq:1})--(\ref{eq:3}) and additional conditions (\ref{eq:4}) 
is well-posed. The corresponding conditions for existence and uniqueness of the solution are available in the above-mentioned works.

In the Hilbert space $H = L_2(\Omega)$, we define the scalar product and norm in the standard way:
\[
  (u,v) = \int_{\Omega} u({\bm x}) v({\bm x}) d{\bm x},
  \quad \|u\| = (u,u)^{1/2} .
\] 
To solve numerically the problem (\ref{eq:1})--(\ref{eq:3}), we employ finite-element 
approximations in space \cite{brenner2008mathematical,Thomee2006}. 
We define the bilinear form
\[
 a(u,v) = \int_{\Omega } \left ( k \, {\rm grad} \, u \, {\rm grad} \, v + c \, u v \right )  d {\bm x} +
 \int_{\partial \Omega } \mu \, u v d {\bm x} .
\] 
We have
\[
 a(u,u) \geq \delta \|u\|^2 , 
 \quad \delta  > 0 .  
\]
Define a subspace of finite elements $V^h \subset H^1(\Omega)$.
Let $\bm x_i, \ i = 1,2, ..., M_h$ be triangulation points for the domain $\Omega$.
For example, when using Lagrange finite elements of the first order (piece-wise linear approximation) we can define pyramid function $\chi_i(\bm x) \subset V^h, \ i = 1,2, ..., M_h$, where
\[
 \chi_i(\bm x_j) = \left \{
 \begin{array}{ll} 
 1, & \mathrm{if~}  i = j, \\
 0, & \mathrm{if~}  i \neq  j .
 \end{array}
 \right . 
\] 
For $v \in V_h$, we have
\[
 v(\bm x) = \sum_{i=i}^{M_h} v_i \chi_i(\bm x),
\] 
where $v_i = v(\bm x_i), \ i = 1,2, ..., M_h$.

We define the discrete elliptic operator $A$ as
\[
 (A y, v) = a(y,v),
\quad \forall \ y,v \in V^h . 
\]
The operator $A$ acts on a finite dimensional space $V^h$ and
\begin{equation}\label{eq:5}
A = A^* \geq \delta I ,
\quad \delta > 0 , 
\end{equation} 
where $I$ is the identity operator in $V^h$.

For the problem (\ref{eq:1})--(\ref{eq:3}), we put into the correspondence the operator equation 
for $w(t) \in V^h$:
\begin{equation}\label{eq:6}
 \frac{d w}{d t} + A w = \varphi, 
 \quad 0 < t \leq T, 
\end{equation} 
\begin{equation}\label{eq:7}
 w(0) =\phi , 
\end{equation} 
where $\varphi  = P f$, $\phi  = P u_0$ with $P$ denoting $L_2$-projection onto $V^h$.
When considering the inverse problem taking into account   (\ref{eq:4}) assume 
\begin{equation}\label{eq:8}
 w(T) = \psi  , 
\end{equation} 
where $\psi   = P u_T$.

\section{Iterative solution of time derivative problem}

For the numerical solution of the inverse problem (\ref{eq:6})--(\ref{eq:8}) with finding $w(t), \varphi$ the simplest approach is to eliminate variable $\varphi$ \cite{samarskii2007numerical,xiangtuan2011direct}.
Differentiating equation (\ref{eq:6}) on time, we obtain
\begin{equation}\label{eq:9}
 \frac{d^2 w}{d t^2} + A \frac{d w}{d t} = 0, 
 \quad 0 < t \leq T . 
\end{equation} 

Further, we consider the boundary value problem (\ref{eq:7})--(\ref{eq:9}).
The correctness of such problem, the computational algorithm and examples of the numerical solution are presented in \cite{samarskii2007numerical}.
The weakness of such approach is caused by the computational complexity of the numerical solution of the boundary problem (\ref{eq:7})--(\ref{eq:9}).
We practically lose the evolutionary character of original problem and must store data in each time step.

The second approach (see, for example, \cite{xiangtuan2011direct}) is based on considering the time derivative.
Let $v = {\displaystyle  \frac{d w}{d t}}$, then equation (\ref{eq:9}) can be written as
\begin{equation}\label{eq:10}
 \frac{d v}{d t} + A v = 0, 
 \quad 0 < t \leq T .
\end{equation} 
For (\ref{eq:10}) we formulate non-local boundary conditions. 
From (\ref{eq:6}) we have
\[
\begin{split}
 v(0) & + A w(0) = \varphi , \\
 v(T) & + A w(T) = \varphi  .
\end{split} 
\] 
Taking into account  (\ref{eq:7}), (\ref{eq:8}) yields
\begin{equation}\label{eq:11}
 v(T) - v(0) = \chi,
 \quad \chi = A(\phi - \psi) .
\end{equation} 

For numerical solution of the problem (\ref{eq:10}), (\ref{eq:11}) we use the simplest iterative refinement of the initial condition for equation (\ref{eq:10}).
The iterative process is organized as follows. The new approximation $k+1$ is found by solving the Cauchy problem: 
\begin{equation}\label{eq:12}
 v^{k+1}(0) =  v^{k}(T) - \chi,
\end{equation} 
\begin{equation}\label{eq:13}
 \frac{d v^{k+1}}{d t} + A v^{k+1} = 0, 
 \quad 0 < t \leq T ,
 \quad k = 0,1, ..., 
\end{equation} 
with some given $v^{0}(0)$.
The desired right-hand side of equation (\ref{eq:6}) is determined using $v^{k+1}(0)$, for example, from the equality
\begin{equation}\label{eq:14}
 \varphi^{k+1} = \phi + v^{k+1}(0) . 
\end{equation} 

When studying the convergence of the iterative process (\ref{eq:12}), (\ref{eq:13})
we consider the problem for error $z^{k+1}(t) = v^{k+1}(t) - v(t)$:
\begin{equation}\label{eq:15}
 z^{k+1}(0) =  z^{k}(T),
\end{equation} 
\begin{equation}\label{eq:16}
 \frac{d z^{k+1}}{d t} + A z^{k+1} = 0, 
 \quad 0 < t \leq T ,
 \quad k = 0,1, ..., 
\end{equation} 
with given $z^{0}(0)$.
Multiplying equation (\ref{eq:16}) for $z^{k}$ in $V^h$ by $z^{k}$, we obtain
\[
 \left ( \frac{d z^{k}}{d t}, z^{k} \right ) + (A z^{k}, z^{k}) = 0 .
\]  
Taking into account (\ref{eq:5}) and  
\[
 \left ( \frac{d z^{k}}{d t}, z^{k} \right ) = \|z^{k}\| \frac{d }{d t} \|z^{k}\|,
\] 
yields
\[
 \frac{d }{d t} \|z^{k}\| + \delta \|z^{k}\| \leq 0 .
\] 
Thus,
\begin{equation}\label{eq:17}
 \|z^{k}(t)\| \leq \exp(-\delta t) \|z^{k}(0)\| .
\end{equation} 
From (\ref{eq:15}) taking into account (\ref{eq:17}) we have
\[
 \|z^{k+1}(0)\| =  \|z^{k}(T)\| \leq \exp(-\delta T) \|z^{k}(0)\| .
\] 
This gives the desired estimate
\begin{equation}\label{eq:18}
 \|v^{k+1}(0) - v(0) \| \leq \varrho \, \|v^{k}(0) - v(0)\|,
 \quad \varrho = \exp(-\delta T),
\end{equation} 
for the convergence of the iterative process (\ref{eq:12})--(\ref{eq:14})  with linear speed $\varrho < 1$.
For the right-hand side with (\ref{eq:14}) we have
\[
 \|\varphi^{k+1} - \varphi\| \leq \varrho \, \|v^{k}(0) - v(0)\| . 
\] 

For computational implementation of proposed algorithm the time approximation deserves special attention.
Let us define a uniform grid in time  
\[
  t_n=n\tau,
  \quad n=0,1,...,N,
  \quad \tau N=T
\]
and denote $y_n = y(t_n), \ t_n = n \tau$.

For the numerical solution of the problem  (\ref{eq:10}), (\ref{eq:11}) we used fully implicit two-level scheme \cite{Samarskiibook,SamarskiiMatusVabischevich2002}, when
\begin{equation}\label{eq:19}
 \frac{v_{n+1} - v_n}{\tau } + A v_{n+1} = 0, 
 \quad n=0,1,...,N-1,
\end{equation} 
\begin{equation}\label{eq:20}
 v_N - v_0 = \chi .
\end{equation} 
The grid problem   (\ref{eq:19}), (\ref{eq:20}) is solved using the following iterative process:
\begin{equation}\label{eq:21}
 v_0^{k+1} = v_N^{k} - \chi ,
\end{equation} 
\begin{equation}\label{eq:22}
 \frac{v_{n+1}^{k+1} - v_n^{k+1}}{\tau } + A v_{n+1}^{k+1} = 0, 
 \quad n=0,1,...,N-1,
 \quad k = 0,1, ..., 
\end{equation} 
where
\begin{equation}\label{eq:23}
 \varphi^{k+1} = \phi + v^{k+1}_0 . 
\end{equation} 

The study of the iterative process (\ref{eq:21})--(\ref{eq:23}) is conducted using the same approach as for the iterative process (\ref{eq:12})--(\ref{eq:14}).
Let now $z_{n+1}^{k+1} = v_{n+1}^{k+1} - v_{n+1}$, then
\begin{equation}\label{eq:24}
 z_0^{k+1} = z_N^{k} ,
\end{equation} 
\begin{equation}\label{eq:25}
 \frac{z_{n+1}^{k+1} - z_n^{k+1}}{\tau } + A z_{n+1}^{k+1} = 0, 
 \quad n=0,1,...,N-1,
 \quad k = 0,1, ... \, . 
\end{equation} 
The key moment of our consideration consists in finding an estimate as (\ref{eq:17}).

We multiply equation (\ref{eq:25}) for $z_{n+1}^{k}$ in $V^h$ by $\tau z_{n+1}^{k}$ and obtain
\[
 \|z_{n+1}^{k}\|^2 + \tau (A z_{n+1}^{k}, z_{n+1}^{k}) = (z_{n}^{k}, z_{n+1}^{k}) .
\] 
Taking into account (\ref{eq:5}) and
\[
 (z_{n}^{k}, z_{n+1}^{k}) \leq \|z_{n}^{k}\| \|z_{n+1}^{k}\|
\] 
we have
\[
 (1 + \tau \delta ) \|z_{n+1}^{k}\| \leq \|z_{n}^{k}\| ,
 \quad n=0,1,...,N-1 .
\] 
Thus,
\begin{equation}\label{eq:26}
 \|z_{n}^{k}\| \leq (1 + \tau \delta )^{-n} \|z_{0}^{k}\| ,
 \quad n=1,2,...,N.
\end{equation} 
Using (\ref{eq:24}), a priori estimate (\ref{eq:26}) allows us to obtain
\[
  \|z_0^{k+1}\| = \|z_N^{k}\| \leq  (1 + \tau \delta )^{-N}  \|z_0^{k}\| .
\] 
Thereby
\begin{equation}\label{eq:27}
 \|v_{0}^{k+1} - v_{0}\| \leq \bar{\varrho} \, \|v_{0}^{k} - v_{0}\|,
 \quad \bar{\varrho} =  (1 + \tau \delta )^{-N} ,
\end{equation} 
which provides the convergence of the iterative process (\ref{eq:21}), (\ref{eq:22}) ($\bar{\varrho} < 1$).
The convergence of the right-hand side with (\ref{eq:23}) is ensured by the estimate
\begin{equation}\label{eq:28}
 \|\varphi^{k+1} - \varphi\| \leq \bar{\varrho}  \, \|v_{0}^{k} - v_{0}\| .
\end{equation} 
This allow us to formulate the following main assertion.

\begin{thm}\label{thm:1}
The iterative process (\ref{eq:21})--(\ref{eq:23}) for the numerical solution of the problem (\ref{eq:6})--(\ref{eq:8}) converges linearly with speed $\bar{\varrho} < 1$ and the estimates (\ref{eq:27}), (\ref{eq:28}) are valid. 
\end{thm}

The outline of the iterative process (\ref{eq:21})--(\ref{eq:23}) is as follows:
\begin{enumerate}
 \item The initial condition is determined using  (\ref{eq:21}). When $k=0$, we assume, for example, $v_0^{0} = - \chi$;
 \item With this initial condition we solve the grid Cauchy problem (\ref{eq:22}). After that we refine the initial condition again. 
\end{enumerate} 

Thus, the computational implementation is based on solving the standard Cauchy problem for parabolic equation on each iteration.

\section{Iterative process for identifying the right-hand side} 

When studying the correctness of the inverse problem (\ref{eq:1})--(\ref{eq:4}) the constructive method of iterative refinement of the right-hand side is often used (see, for example, \cite{prilepko1993certain,PrilepkoKostin,Kostin2013}). We consider the possibility of using this approach for the approximate solution of the problem (\ref{eq:6})--(\ref{eq:8}).

In the new iterative step the right-hand side is determined using (\ref{eq:6}) for $t = T$, taking into account (\ref{eq:8}):
\begin{equation}\label{eq:29}
 \varphi^{k+1} = \frac{d w^k}{d t}(T) + A \psi , 
 \quad k = 0,1, ... ,
\end{equation} 
with some given initial assumption $\varphi^{0}$.
Then, the Cauchy problem is solved:
\begin{equation}\label{eq:30}
 \frac{d w^{k+1}}{d t} + A w^{k+1} = \varphi^{k+1}, 
 \quad 0 < t \leq T, 
\end{equation} 
\begin{equation}\label{eq:31}
 w^{k+1}(0) =\phi .
\end{equation} 

To estimate the convergence we introduce
\[
 \eta^{k+1} = \varphi^{k+1} - \varphi ,
 \quad z^{k+1}(t) = w^{k+1}(t) - w(t),
 \quad 0 \leq t \leq T, 
 \quad k = 0,1, ... \, .
\] 
Then, from (\ref{eq:29})--(\ref{eq:31}) we obtain
\begin{equation}\label{eq:32}
 \eta^{k+1} = \frac{d z^k}{d t}(T) , 
 \quad k = 0,1, ... ,
\end{equation} 
\begin{equation}\label{eq:33}
 \frac{d z^{k+1}}{d t} + A z^{k+1} = \eta^{k+1}, 
 \quad 0 < t \leq T, 
\end{equation} 
\begin{equation}\label{eq:34}
 z^{k+1}(0) = 0 .
\end{equation} 
We differentiate  equation (\ref{eq:33}) with respect to time and rewrite it for 
${\displaystyle v^k = \frac{d z^{k}}{d t}}$ in the form
\begin{equation}\label{eq:35}
 \frac{d v^{k}}{d t} + A v^{k} = 0, 
 \quad 0 < t \leq T, 
\end{equation} 
The initial condition for  (\ref{eq:35}) we obtain from (\ref{eq:33}) with
 $t= 0$, using (\ref{eq:34}):
\begin{equation}\label{eq:36}
 v^k(0) =  \eta^{k} .
\end{equation} 
Similar to (\ref{eq:17}) we have the estimate
\begin{equation}\label{eq:37}
 \|v^k(t)\| \leq \exp(-\delta t) \|\eta^{k}\|
\end{equation} 
for the solution of the problem (\ref{eq:35}), (\ref{eq:36}).
In view of (\ref{eq:37}) from (\ref{eq:32}) we have the estimate
\[
 \|\eta^{k+1}\| = \|v^k(T)\| \leq \exp(-\delta T) \|\eta^{k}\|, 
 \quad k = 0,1, ... 
\] 
for the convergence of desired right-hand side.

The time discretization is again formulated in the basis of implicit approximation.  
Formally, we define the solution of the problem (\ref{eq:33}), (\ref{eq:34}) on expanded grid:
\[
  t_n=n\tau,
  \quad n=-1,0,...,N,
  \quad \tau N=T .
\]
From (\ref{eq:6})--(\ref{eq:8}) we come to the problem
\begin{equation}\label{eq:38}
 \frac{w_{n+1} - w_{n}}{\tau } + A w_{n+1} = \varphi, 
 \quad n=-1,0,...,N-1 
\end{equation} 
\begin{equation}\label{eq:39}
 w_0 = \phi ,
\end{equation} 
\begin{equation}\label{eq:40}
 w_N = \psi .
\end{equation} 

For the numerical solution of the problem (\ref{eq:38})--(\ref{eq:40}) the grid analogue of the iterative process (\ref{eq:29})--(\ref{eq:31}) is applied.
Now, we compare the relationship
\begin{equation}\label{eq:41}
 \varphi^{k+1} = \frac{w^k_N - w^k_{N-1}}{\tau} + A \psi , 
 \quad k = 0,1, ... , 
\end{equation} 
with (\ref{eq:29}). 
To approximate equation (\ref{eq:30}) the implicit difference scheme is used 
\begin{equation}\label{eq:42}
 \frac{w^{k+1}_{n+1} - w^{k+1}_{n}}{\tau } + A w^{k+1}_{n+1} = \varphi^{k+1}, 
 \quad n=-1,0,...,N-1 ,
\end{equation} 
under condition  (see  (\ref{eq:31}))
\begin{equation}\label{eq:43}
 w^{k+1}_0 =\phi .
\end{equation} 

\begin{thm}\label{thm:2}
The iterative process (\ref{eq:41})--(\ref{eq:43}) for the numerical solution of the problem (\ref{eq:38})--(\ref{eq:40}) converges linearly with speed 
\[
 \bar{\varrho} = (1 + \tau \delta )^{-N} ,
\] 
and the estimate
\begin{equation}\label{eq:44}
 \|\varphi^{k+1} - \varphi\| \leq \bar{\varrho}  \,  \|\varphi^{k} - \varphi\|.
\end{equation} 
is valid.
\end{thm}

\begin{pf}
We define
\[
 \eta^{k+1} = \varphi^{k+1} - \varphi ,
 \quad z^{k+1}_{n+1} = w^{k+1}_{n+1} - w_{n+1},
 \quad n = -1,0, ..., N-1,
 \quad k = 0,1, ... \, .
\] 
From (\ref{eq:41})--(\ref{eq:43}) we obtain
\begin{equation}\label{eq:45}
 \eta^{k+1} = \frac{z^k_N - z^k_{N-1}}{\tau} , 
 \quad k = 0,1, ... ,
\end{equation} 
\begin{equation}\label{eq:46}
 \frac{z^{k+1}_{n+1} - z^{k+1}_{n}}{\tau } + A z^{k+1}_{n+1} = \eta^{k+1}, 
 \quad n=-1,0,...,N-1 ,
\end{equation} 
\begin{equation}\label{eq:47}
 z^{k+1}_0 = 0 .
\end{equation} 
Let's assume
\[
 v^k_{n+1} = \frac{z^{k}_{n+1} - z^{k}_{n}}{\tau } ,
 \quad n=-1,0,...,N-1 .
\] 
From (\ref{eq:46}) we get
\begin{equation}\label{eq:48}
 \frac{v^{k}_{n+1} - v^{k}_{n}}{\tau } + A v^{k}_{n+1} = 0, 
 \quad n= 0,1,...,N-1 .
\end{equation} 
From equation (\ref{eq:46}) when $n=0$ using the condition (\ref{eq:47}) we have
\begin{equation}\label{eq:49}
 v^{k}_{0} = \eta^k . 
\end{equation} 
Similar to (\ref{eq:26}) for the solution of the problem (\ref{eq:48}), (\ref{eq:49}) we prove the a priori estimate
\begin{equation}\label{eq:50}
 \|v_{n}^{k}\| \leq (1 + \tau \delta )^{-n} \|\eta^{k}\| ,
 \quad n=1,2,...,N.
\end{equation} 
From (\ref{eq:45}) and (\ref{eq:50}) taking into account introduced notation  yields
\[
 \|\eta^{k+1}\| = \|v_{N}^{k}\|  \leq (1 + \tau \delta )^{-N} \|\eta^{k}\| , 
 \quad k = 0,1, ... \, .
\] 
Thus, we establish the estimate (\ref{eq:44}).
\end{pf}

\section{Integral overdetermination} 

When considering inverse problem of identifying the right-hand side of parabolic equation, the integral overdetermination is often used instead of the final  overdetermination.
In this case, in place of equation (\ref{eq:4}) the following condition is involved
\begin{equation}\label{eq:51}
  \int_{0}^{T}\omega(t) u({\bm x}, t) d t = u_{T} ({\bm x}),
  \quad {\bm x} \in \Omega ,
\end{equation} 
where $\omega(t)$ --- given function and
\[
 \omega(t) \geq 0,
 \quad  \int_{0}^{T}\omega(t) dt = 1.
\] 
For the numerical solution of the inverse problem (\ref{eq:1})--(\ref{eq:3}), (\ref{eq:51})
iterative process considered above can be used. 
Here, we note capabilities of iterative refinement of the desired right-hand side 
\cite{prilepko1993certain,PrilepkoKostin}.

After finite-element approximation in space we come to equation  (\ref{eq:6}),
which is supplemented with the initial condition (\ref{eq:7}). Taking into account (\ref{eq:51}), we have  
\begin{equation}\label{eq:52}
  \int_{0}^{T}\omega(t) w(t) d t = \psi .
\end{equation} 
instead of (\ref{eq:8}).
Integrating equation  (\ref{eq:6}) with weight $\omega(t)$ over $t$ from $0$ to $T$ and taking into account (\ref{eq:52}),
we obtain
\begin{equation}\label{eq:53}
  \int_{0}^{T}\omega(t)\frac{d w}{d t}(t) d t + A \psi = \varphi .
\end{equation} 

The iterative refinement of the right-hand side is based on the relation  (\ref{eq:53}) (see (\ref{eq:29})):
\begin{equation}\label{eq:54}
 \varphi^{k+1} = \int_{0}^{T}\omega(t) \frac{d w^k}{d t}(t) d t + A \psi , 
 \quad k = 0,1, ... \, .
\end{equation} 
With new approximation for the right-hand side the Cauchy problem  (\ref{eq:30}),  (\ref{eq:51}) is solved.
The study of the convergence of the iterative process is performed similarly to the study of the convergence of the process (\ref{eq:29}), (\ref{eq:51}).

In view of early introduced notation we have
\begin{equation}\label{eq:55}
 \eta^{k+1} = \int_{0}^{T}\omega(t) \frac{d z^k}{d t}(t) d t , 
 \quad k = 0,1, ... ,
\end{equation} 
and for $z^{k+1}$ we have the problem (\ref{eq:33}),  (\ref{eq:34}).
For ${\displaystyle v^k = \frac{d z^{k}}{d t}}$ we have the problem (\ref{eq:35}),  (\ref{eq:36})
and the estimate (\ref{eq:37}) exists.
Taking into account this from  (\ref{eq:5}) we obtain
\[
\begin{split}
 \|\eta^{k+1}\| & = \left  \|\int_{0}^{T}\omega(t) v^k(t) d t \right  \| =
 \int_{0}^{T}\omega(t) \|v^k(t)\| d t  \\
 & \leq  \int_{0}^{T}\omega(t) \exp(-\delta t) d t \  \|\eta^{k}\| ,
 \quad k = 0,1, ... \, .
\end{split}
\] 
Thus,
\begin{equation}\label{eq:56}
 \|\varphi^{k+1} - \varphi\| \leq \varrho \, \|\varphi^{k} - \varphi \| ,
 \quad k = 0,1, ... ,
\end{equation} 
at that
\[
  \varrho = \int_{0}^{T}\omega(t) \exp(-\delta t) d t .
\] 
For usual assumptions about continuity of the weight function  $\omega(t)$ 
we have $\varrho < 1$.
In the limit case $\omega(t) = \delta(t-T)$  ($\delta(t)$ --- $\delta$-function)
from (\ref{eq:52}) we obtain the condition of final overdetermination (\ref{eq:8}).
Other limit case, when $\omega(t) = \delta(t)$, is not interesting, since two condition for $t=0$ are set and, as expected, the iterative process does not converge ($\varrho = 1$).

\begin{thm}\label{thm:3}
The iterative process (\ref{eq:30}), (\ref{eq:31}), (\ref{eq:55}) for the numerical solution of the problem (\ref{eq:6}), (\ref{eq:7}), (\ref{eq:52}) for (\ref{eq:5}) converges linearly with the speed 
$\varrho$, while the estimate (\ref{eq:56}) is correct.
\end{thm}

Applying  an approach similar to one in Section 3 the iterative process is studied when implicit time discretization is used.  

\section{More general problems} 

We have investigated the iterative methods for the numerical solution of the inverse problem (\ref{eq:1})--(\ref{eq:4}), when the right-hand side does not depend on time. 
In more general case, the problem of identifying multiplicative right-hand side, when the dependence of the right-hand side on time is known and the dependence on spatial variables is unknown, is stated.

Instead of  (\ref{eq:1}) we consider the following equation 
\begin{equation}\label{eq:57}
  \frac{\partial u}{\partial t}- {\rm div} (k({\bm x}) {\rm grad}u) + c({\bm x}) u 
  = \beta(t) f({\bm x}),
  \quad {\bm x} \in \Omega,
  \quad 0 < t \leq T ,   
\end{equation} 
where $\beta(t)$ --- some given function.
The inverse problem of finding the pair $u(\bm x,t), \ f(\bm x)$ 
from (\ref{eq:2})--(\ref{eq:4}), (\ref{eq:57}) is stated.
We assume that
\begin{equation}\label{eq:58}
 \beta(t) > 0,
 \quad \frac{d \beta }{d t} \geq 0,
 \quad 0 \leq t \leq T,
 \quad \beta(T) = 1 .  
\end{equation} 

After approximation in space we come to the equation
\begin{equation}\label{eq:59}
 \frac{d w}{d t} + A w = \beta(t) \varphi, 
 \quad 0 < t \leq T, 
\end{equation} 
for $w(t) \in V^h$, which is supplemented by (\ref{eq:7}), (\ref{eq:8}).
Taking into account condition (\ref{eq:8}) and normalization $\beta(T) = 1$ equation (\ref{eq:59})
for $t=T$ gives
\begin{equation}\label{eq:60}
 \varphi = \frac{d w}{d t} (T) + A \psi .
\end{equation} 

Iterative refinement of the right-hand side is performed using (\ref{eq:60}) and has the form
(\ref{eq:29}). When  right-hand side is known on the new iteration we solve equation
\begin{equation}\label{eq:61}
 \frac{d w^{k+1}}{d t} + A w^{k+1} = \beta(t) \varphi^{k+1}, 
 \quad 0 < t \leq T, 
\end{equation} 
with the initial condition (\ref{eq:31}).
For errors $\eta^{k+1}$ and $z^{k+1}(t)$  
we obtain the relation (\ref{eq:29}), and equation
\begin{equation}\label{eq:62}
 \frac{d z^{k+1}}{d t} + A z^{k+1} = \beta (t) \eta^{k+1}, 
 \quad 0 < t \leq T, 
\end{equation} 
with a uniform initial condition  (\ref{eq:34}).

As before when considering the convergence of the iterative process (\ref{eq:29})--(\ref{eq:31}),
we study the problem for ${\displaystyle v^k = \frac{d z^{k}}{d t}}$. From  (\ref{eq:62})
we have
\begin{equation}\label{eq:63}
 \frac{d v^{k}}{d t} + A v^{k} = \frac{d \beta}{d t}(t) \eta^{k+1}, 
 \quad 0 < t \leq T .
\end{equation} 
From equation (\ref{eq:62}) for $t= 0$ taking into account (\ref{eq:34}) the initial condition for equation (\ref{eq:63}) is obtained:
\begin{equation}\label{eq:64}
 v^k(0) =  \beta (0) \eta^{k} .
\end{equation} 

We multiply equation (\ref{eq:63}) by $V^h$ on $v^{k}$, which leads to
\[
 \left  ( \frac{d v^{k}}{d t}, v^{k} \right )  + 
 (A v^{k}, v^{k}) = \frac{d \beta}{d t}(t) (\eta^{k}, v^{k}) .  
\] 
Taking into account that the derivative of the function $\beta(t)$ is non negative and inequality 
\[
 (\eta^{k}, v^{k}) \leq  \|\eta^{k}\| \, \|v^{k}\|,
\]  
we obtain
\[
 \frac{d }{d t} \|v^{k}\| + \delta \|v^{k}\| \leq  \frac{d \beta}{d t}(t) \|\eta^{k}\|.
\] 
This implies
\[
 \frac{d }{d t} (\exp (\delta t) \|v^{k}\|) \leq  
 \exp (\delta t) \frac{d \beta}{d t}(t) \|\eta^{k}\|.
\] 
Integration by $t$ from $0$ to $T$ yields
\[
 \exp (\delta T) \|v^{k}(T)\| - \|v^{k}(0)\| \leq 
 \int_{0}^{T}\exp (\delta t) \frac{d \beta}{d t}(t) d t \ \|\eta^{k}\| .
\] 
Given the initial condition (\ref{eq:64}) and
\[
 \int_{0}^{T}\exp (\delta (t-T)) \frac{d \beta}{d t}(t) d t \leq 
 \int_{0}^{T} \frac{d \beta}{d t}(t) d t = 1 - \beta(0), 
\] 
we have
\begin{equation}\label{eq:65}
 \|z^{k}(T)\| \leq \varrho  \|z^{k}(0)\| ,
\end{equation} 
where
\begin{equation}\label{eq:66}
 \varrho = 1 - \beta(0) (1 - \exp (- \delta T)) .
\end{equation} 
From (\ref{eq:29}) and (\ref{eq:65}) we obtain the estimate (\ref{eq:65}), where
$\varrho < 1$  is determined according to (\ref{eq:66}).
The result of consideration is the following statement.

\begin{thm}\label{thm:4}
The iterative process (\ref{eq:29}), (\ref{eq:30}), (\ref{eq:61}) for the numerical solution of the problem (\ref{eq:7}), (\ref{eq:8}), (\ref{eq:59}) for (\ref{eq:5}) converges linearly at a rate of $\varrho$, which is determined according to (\ref{eq:66}), while the estimate (\ref{eq:65}) is correct.
\end{thm}

\section{Numerical experiments} 

\begin{figure}[htp]
  \begin{center}
    \includegraphics[scale = 0.5] {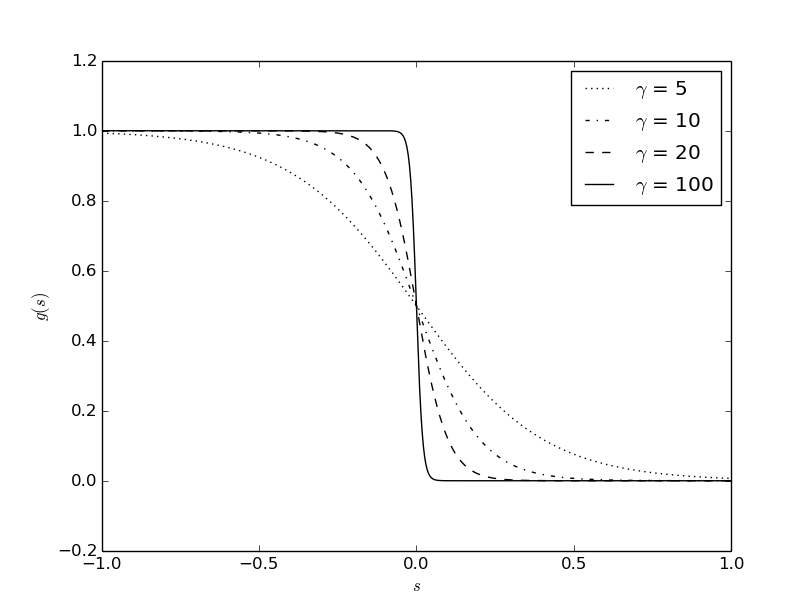}
	\caption{Function $g(s) = (1 + \exp(\gamma s))^{-1}$ at different values $\gamma$}
	\label{fig:1}
  \end{center}
\end{figure}

We illustrate the capabilities of iterative solution of inverse problems of identification of the right-hand side of parabolic equations by results of the numerical solution of a test problem.
We consider model problem, when
\begin{equation}\label{eq:67}
  \frac{\partial u}{\partial t}- {\rm div} \, {\rm grad} \, u + c u = f({\bm x}),
  \quad {\bm x} \in \Omega,
  \quad 0 < t \leq T , 
\end{equation}   
\begin{equation}\label{eq:68}
  \frac{\partial u}{\partial n} = 0,
  \quad {\bm x} \in \partial\Omega,
  \quad 0 < t \leq T, 
\end{equation} 
\begin{equation}\label{eq:69}
  u({\bm x}, 0) = 0,
  \quad {\bm x} \in \Omega . 
\end{equation} 
The forward problem is solved in the unit square
\[
 \Omega = \{ \bm x = (x_1, x_2) \ | \ 0 < x_1 < 1, \ 0 < x_2 < 1 \} 
\] 
with given right-hand side $f({\bm x})$ and the solution at the finite time is found 
\[
 u_T(\bm x) = u(\bm x,T) .
\] 
After that, the inverse problem is solved when $u_T(\bm x)$ is known, but we need to find $f(\bm x)$.
The right-hand side is taken as
\begin{equation}\label{eq:70}
 f(\bm x) = \frac{1}{1 + \exp(\gamma (x_1-x_2))} .
\end{equation} 
For $\gamma \rightarrow \infty$, the right-hand side becomes discontinuous (Fig.\ref{fig:1}).

The forward problem (\ref{eq:67})--(\ref{eq:70}) is solved within the first quasi-real computational experiment. The solution of this problem at the finite time (the function $ u(\bm x,T)$) is used as input data for the inverse problem. 
In our analysis we focus on the iterative solution of the problem of identification after the finite element approximation in space.
Because of this, we do not discuss the dependence of the accuracy of the approximate solution on the approximation in the space, which is more appropriate in another study.
We perform the evaluation of the effect of computational errors on the basis of calculations on different time grids, when using the input data derived from the solution of the forward problem on more detailed time grid and with a more accurate approximations in time.

\begin{figure}[htp]
  \begin{center}
    \includegraphics[scale = 0.25] {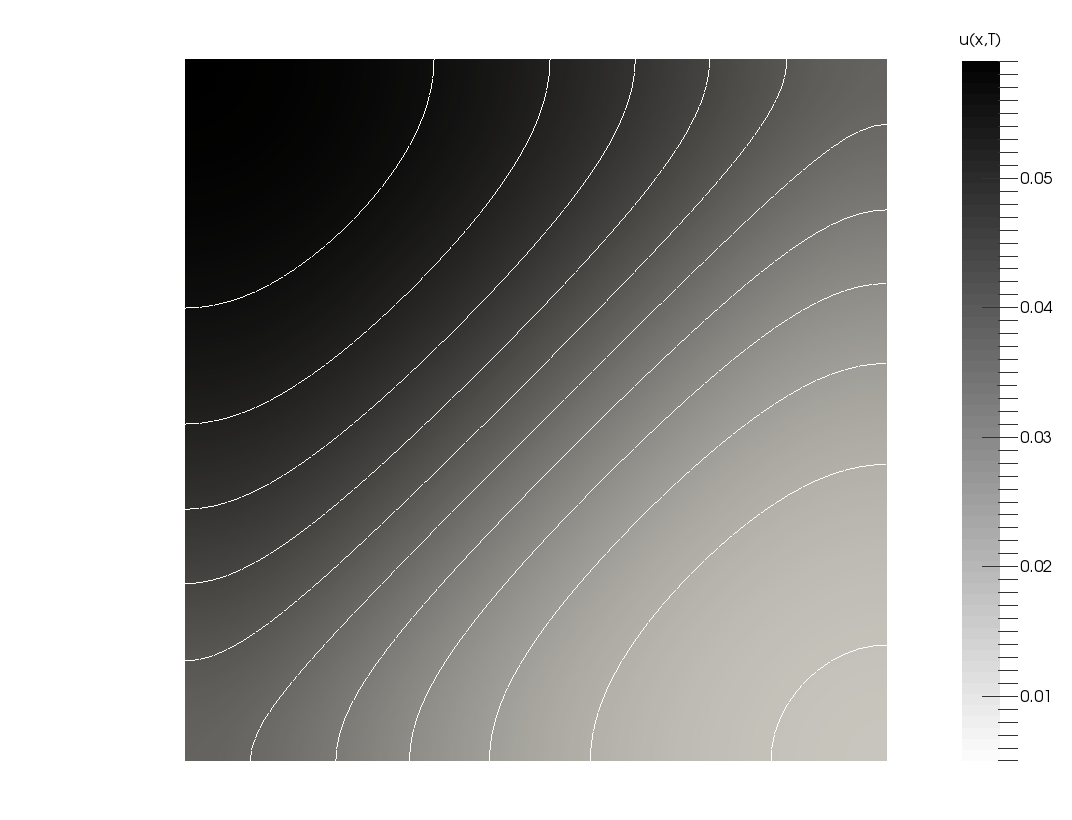}
	\caption{The solution of the forward problem $u_T(\bm x) = u(\bm x,T)$}
	\label{fig:2}
  \end{center}
\end{figure}

For the base case we set $c = 10$, $T = 0.1$, $\gamma =10$.
When solving the forward problem we use the Crank-Nicolson scheme for time discretization, the time step is $\tau = 1\cdot 10^{-4}$.
The uniform mesh with the division into 50 intervals in each direction is used, the Lagrangian finite elements of the second degree are applied.
The solution at the finite time is shown in Fig.\ref{fig:2}.

The inverse problem is solved using fully implicit scheme (see  (\ref{eq:22}), (\ref{eq:30})).
The error of the approximate solution of the problem of identification on a single iteration is evaluated as follows
\[
 \varepsilon_\infty(k) = \max_{\bm x \in \Omega} |\varphi^k(\bm x) - f(\bm x)|,
\] 
\[
 \varepsilon_2(k) = \|\varphi^k(\bm x) - f(\bm x)\|,
\] 
where $\varphi(\bm x)$ --- the approximate solution, and $f(\bm x)$ --- the exact solution of the inverse problem.

\begin{figure}[htp]
  \begin{center}
\begin{minipage}{0.49\linewidth}
\center{\includegraphics[width=1\linewidth]{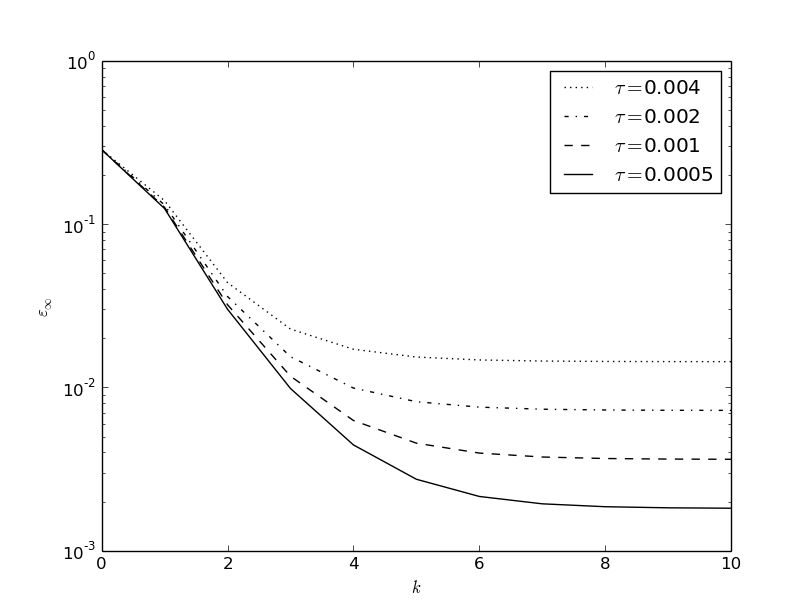}} \\
\end{minipage}
\hfill
\begin{minipage}{0.49\linewidth}
\center{\includegraphics[width=1\linewidth]{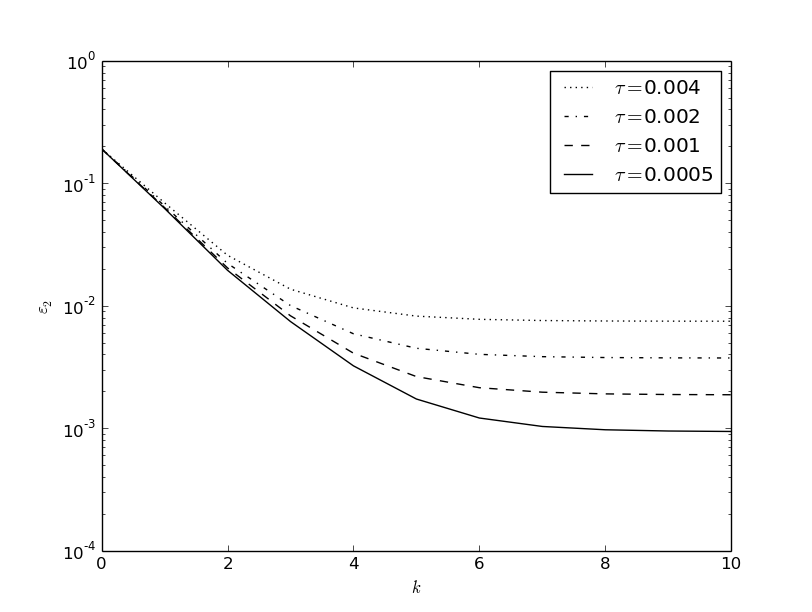}} \\
\end{minipage}
\caption{The iterative process  (\ref{eq:21})--(\ref{eq:23})} 
\label{fig:3}
  \end{center}
\end{figure}

\begin{figure}[htp]
  \begin{center}
\begin{minipage}{0.49\linewidth}
\center{\includegraphics[width=1\linewidth]{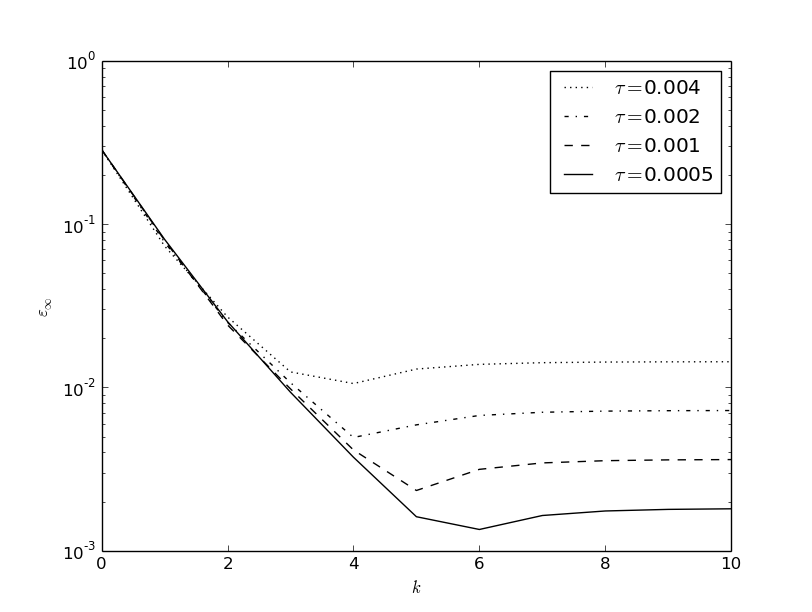}} \\
\end{minipage}
\hfill
\begin{minipage}{0.49\linewidth}
\center{\includegraphics[width=1\linewidth]{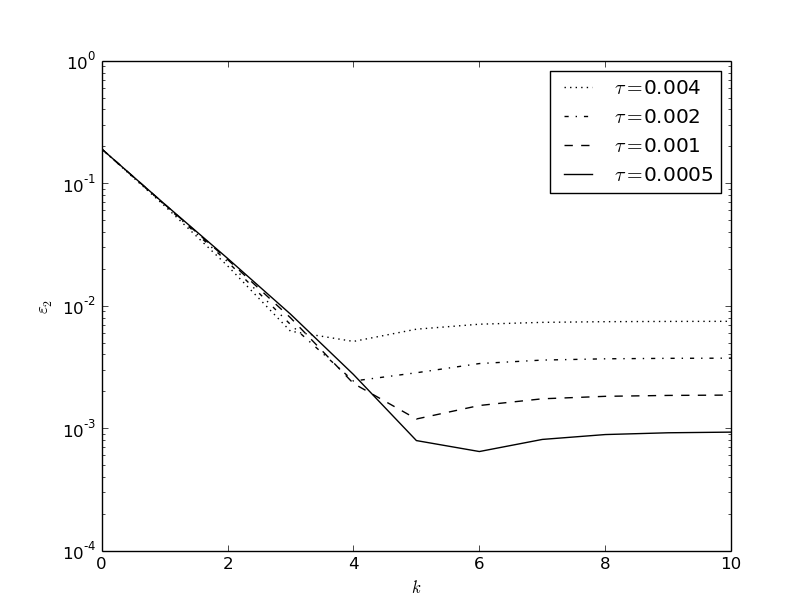}} \\
\end{minipage}
\caption{The iterative process (\ref{eq:41})--(\ref{eq:43})} 
\label{fig:4}
  \end{center}
\end{figure}

Influence of time step of the iterative process (\ref{eq:21})--(\ref{eq:23}) on accuracy is illustrated in Fig.~\ref{fig:3}.
We see the rapid convergence of the iterative process and improving of the accuracy of the approximate solution by reducing the time step. 
Similar data for the iterative process (\ref{eq:41})--(\ref{eq:43}) are shown in Fig.\ref{fig:4}. 
The initial approximation is taken in the form
\[
 v^0(0) = A \psi ,
 \quad \varphi^0 = A \psi . 
\] 

\begin{figure}[htp]
  \begin{center}
\begin{minipage}{0.49\linewidth}
\center{\includegraphics[width=1\linewidth]{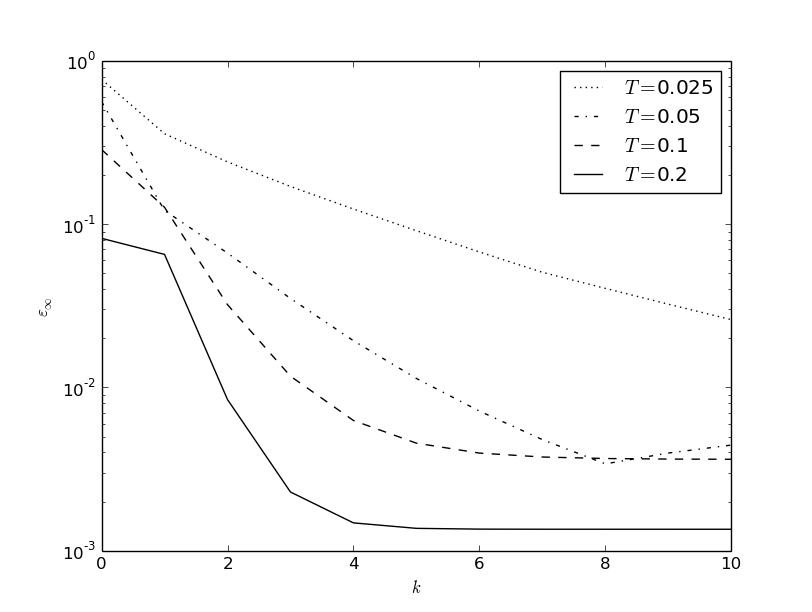}} \\
\end{minipage}
\hfill
\begin{minipage}{0.49\linewidth}
\center{\includegraphics[width=1\linewidth]{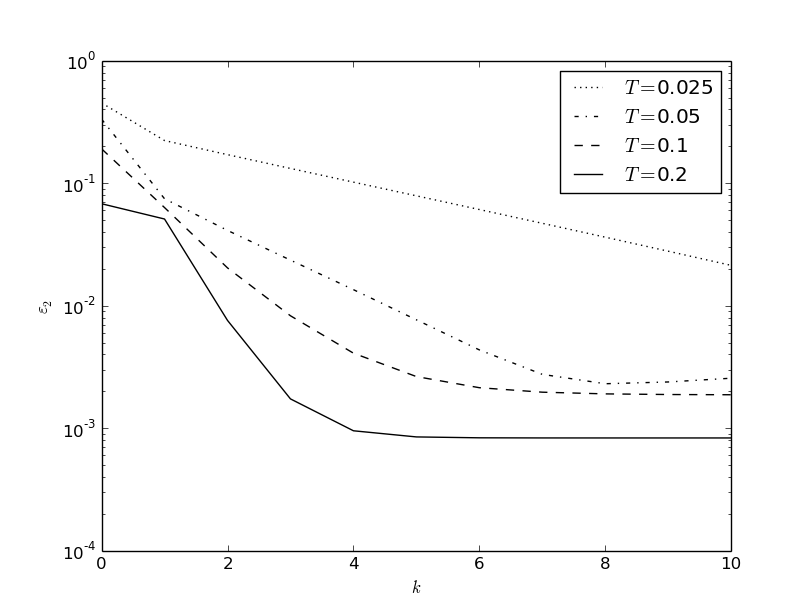}} \\
\end{minipage}
\caption{The dependence on $T$: the iterative process (\ref{eq:21})--(\ref{eq:23})} 
\label{fig:5}
  \end{center}
\end{figure}

The rate of convergence of iterative processes (\ref{eq:21})--(\ref{eq:23}) and (\ref{eq:41})--(\ref{eq:43}) depends essentially on the observation interval --- see assessment (\ref{eq:28}) and (\ref{eq:44}).
The relevant data of numerical computations are presented in Fig. \ref{fig:5} and \ref{fig:6}. 
These computations are carried out at $\tau = 1\cdot 10^{-3}$.

\begin{figure}[htp]
  \begin{center}
\begin{minipage}{0.49\linewidth}
\center{\includegraphics[width=1\linewidth]{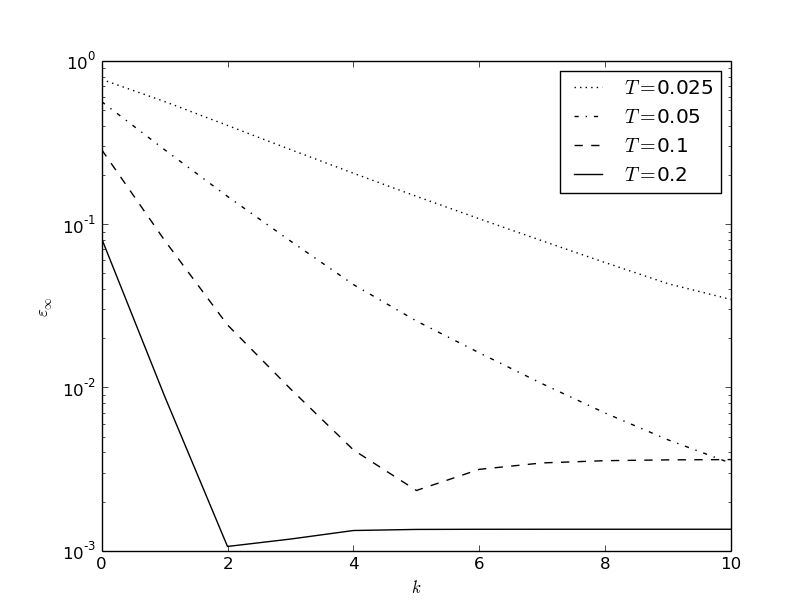}} \\
\end{minipage}
\hfill
\begin{minipage}{0.49\linewidth}
\center{\includegraphics[width=1\linewidth]{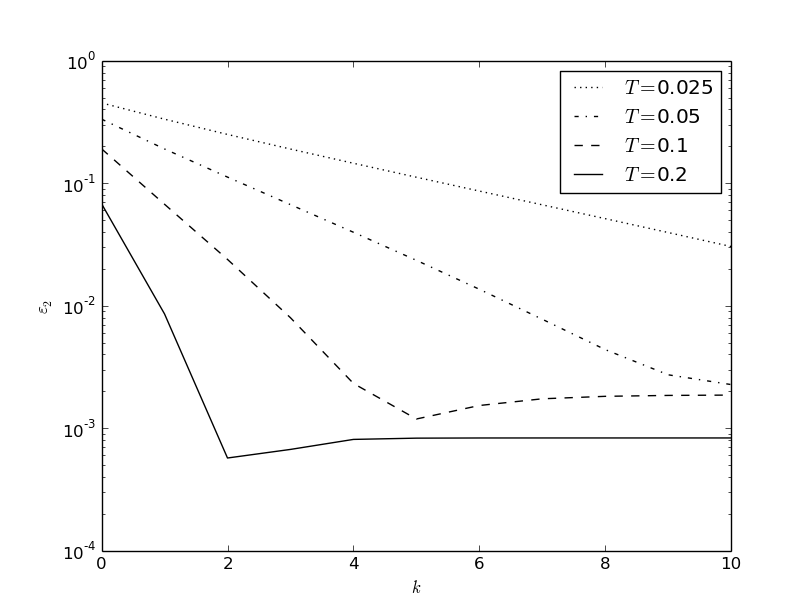}} \\
\end{minipage}
\caption{The dependence on $T$: the iterative process (\ref{eq:41})--(\ref{eq:43})} 
\label{fig:6}
  \end{center}
\end{figure}

\begin{figure}[htp]
  \begin{center}
\begin{minipage}{0.49\linewidth}
\center{\includegraphics[width=1\linewidth]{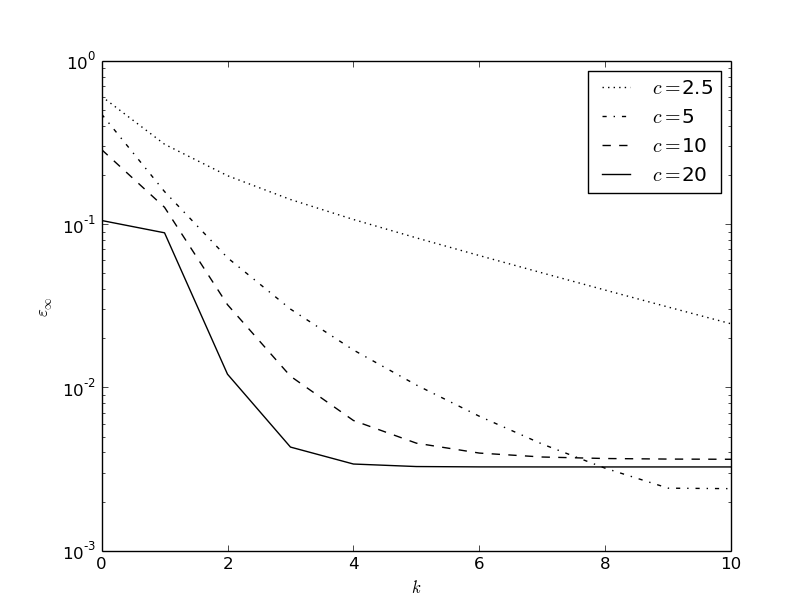}} \\
\end{minipage}
\hfill
\begin{minipage}{0.49\linewidth}
\center{\includegraphics[width=1\linewidth]{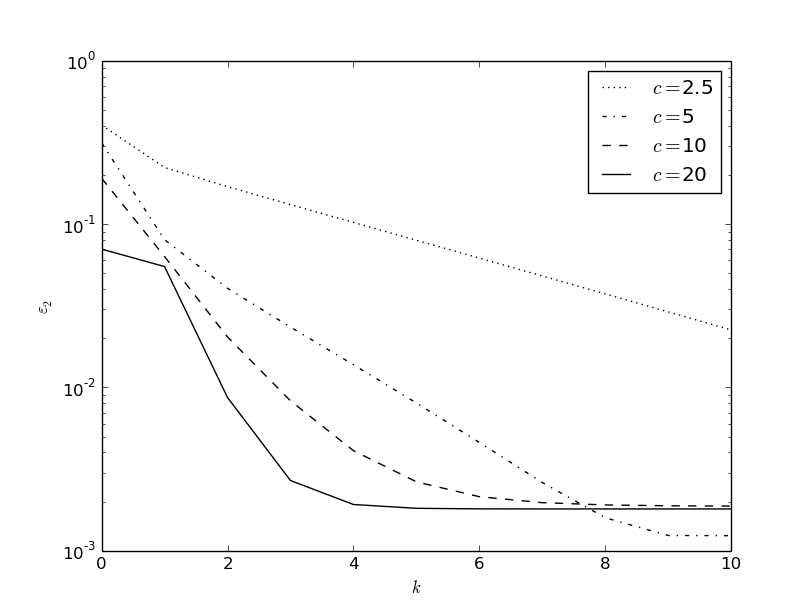}} \\
\end{minipage}
\caption{The dependence on $c$: the iterative process (\ref{eq:41})--(\ref{eq:43})} 
\label{fig:6}  
\label{fig:7}
  \end{center}
\end{figure}

\begin{figure}[htp]
  \begin{center}
\begin{minipage}{0.49\linewidth}
\center{\includegraphics[width=1\linewidth]{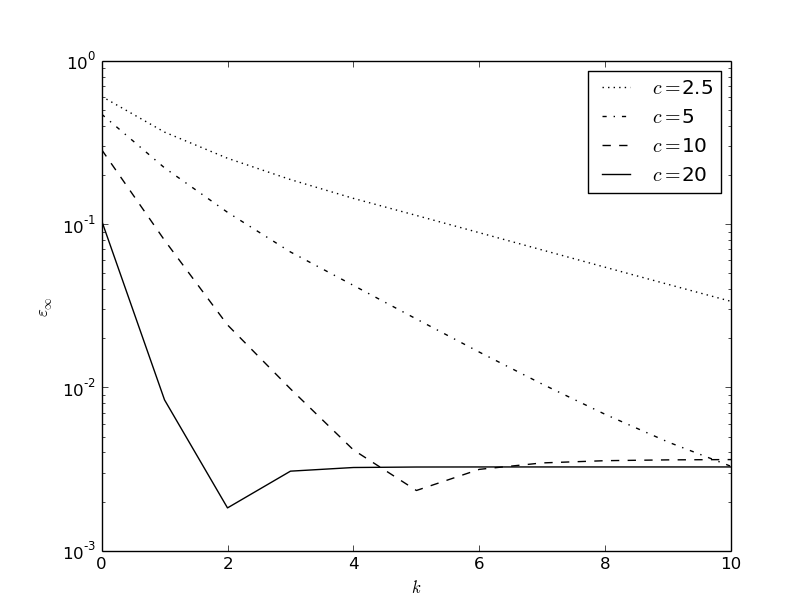}} \\
\end{minipage}
\hfill
\begin{minipage}{0.49\linewidth}
\center{\includegraphics[width=1\linewidth]{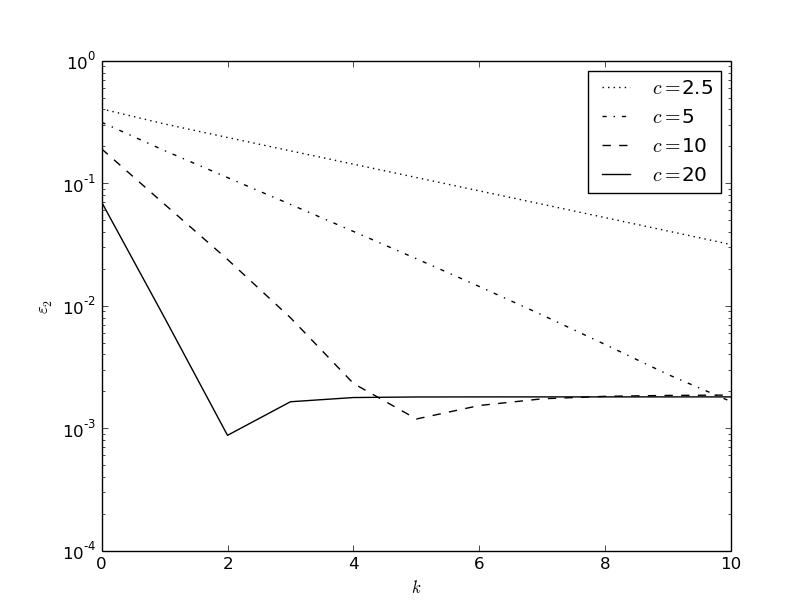}} \\
\end{minipage}
\caption{The dependence on $c$: the iterative process (\ref{eq:41})--(\ref{eq:43})} 
\label{fig:6}  
\label{fig:8}
  \end{center}
\end{figure}

There is a similar dependence of the rate of convergence of the iterative processes (\ref{eq:21})--(\ref{eq:23}) and (\ref{eq:41})--(\ref{eq:43}) on the value of $\delta$. 
For our test problem (\ref{eq:67})--(\ref{eq:69}) we have $\delta = c$.
The convergence for different values $c$ is illustrated in Figs \ref{fig:7} and \ref{fig:8}. 
These computations are carried out at $\tau = 1\cdot 10^{-3}$ for $T = 0.1$.

The accuracy of the approximate solution, the rate of convergence of the iterative processes (\ref{eq:21})--(\ref{eq:23}) and (\ref{eq:41})--(\ref{eq:43})  weakly depends on the required right-hand side.
This is particularly evidenced by the data presented in Tables \ref{tab:1} and \ref{tab:2}.
We consider the problem of identification with the exact solution (\ref{eq:70}) for different values $\gamma$.

\begin{table}
\caption{The dependence of $\varepsilon_\infty$ on $\gamma$: the iterative process (\ref{eq:21})--(\ref{eq:23}).} 
\begin{center}
{\begin{tabular}{@{}ccccc}
  \hline
$k$ & $\gamma = 5$ & $\gamma = 10$ & $\gamma = 20$ & $\gamma = 100$ \\
  \hline
0 & 0.2698616 & 0.2859782  &     0.2915919  &    0.2935690  \\   
1 & 0.1147344 & 0.1266791  &     0.1307726  &    0.1321993  \\   
2 & 0.0300993 & 0.0320972  &     0.0327864  &    0.0330277  \\   
3 & 0.0111847 & 0.0117121  &     0.0118796  &    0.0119350  \\   
4 & 0.0059519 & 0.0062735  &     0.0063679  &    0.0063972  \\   
5 & 0.0042635 & 0.0045563  &     0.0046405  &    0.0046662  \\   
  \hline
\end{tabular}}
\label{tab:1}
\end{center}
\end{table}

\begin{table}
\caption{The dependence of $\varepsilon_2$ on $\gamma$: the iterative process (\ref{eq:21})--(\ref{eq:23}).}
\begin{center}
{\begin{tabular}{@{}ccccc}
\hline
$k$ & $\gamma = 5$ & $\gamma = 10$ & $\gamma = 20$ & $\gamma = 100$ \\
\hline
0 & 0.1889483 & 0.1910198  & 0.1918367  & 0.1921380 \\   
1 & 0.0596180 & 0.0633147  & 0.0647253  & 0.0652396 \\   
2 & 0.0200653 & 0.0203738  & 0.0204951  & 0.0205398 \\   
3 & 0.0082159 & 0.0082824  & 0.0083087  & 0.0083185 \\   
4 & 0.0040403 & 0.0041006  & 0.0041245  & 0.0041333 \\   
5 & 0.0025596 & 0.0026414  & 0.0026734  & 0.0026853 \\   
\hline
\end{tabular}}
\label{tab:2}
\end{center}
\end{table}

\section{Conclusions} 

\begin{enumerate}
 \item The problem of identifying the source of a parabolic equation of the second order, which depends only on the space variables.
Additional information about the solution are set at the final time.
To approximate in space we use the standard Lagrange finite elements for the triangulation of the computational domain.
 \item The original inverse problem is formulated as a problem with a non-local condition for the time for evolution equation of the first order.
 For the approximate solutions of non-classical boundary value problem we propose the iterative process of refinement of the initial condition.
We have investigated the rate of convergence of the iterative process with respect to the permanent positive definiteness of operator of the problem  and  interval time.
 \item We constructed the numerical algorithm based on the well-known approach with iterative refinement of the right-hand side. This algorithm have previously been used to prove the unique solvability of the inverse problem.
We have investigated the rate of convergence of the iterative process for an operator differential equation of first order (semi-discrete approximation).
Similar results are obtained when using a fully implicit approximations of time (fully discrete approximation).
 \item We have considered the inverse problem of identification of the right side of a parabolic equation when more information about the decision is given in the form of a weighted average over time.
The convergence of the iterative process with the refinement of the required right-hand side and the rate of convergence is shown.
 \item We have studied more general problem with multiplicative right-hand side, when one multiplier determines the known dependence of the right-hand side on the time, and the second defines the unknown dependence on the spatial variables.
 The estimates of the rate of convergence of the iterative process with refinement of the right-hand side with respect to parameters of the problem are obtained.
 \item The capabilities of considered computational algorithms of identification of the right-hand side are illustrated by  the results of a numerical solution of the model inverse problem for the two-dimensional parabolic equation in the square.
Within the quasi-real computational experiment the additional data for the inverse problem are found using the results of the numerical solution of the forward problem.
The calculated data show a high rate of convergence of the discussed iterative processes of Picard type which do not contain any adjusting parameters.
To find the approximate solution of the inverse problem it is enough to perform 5, 6 iterations, where the standard boundary value problem for a parabolic equation is solved.  
\end{enumerate} 

\subsection*{Acknowledgements}

This research was supported by RFBR (project 14-01-00785).

\end{document}